\documentclass{amsart}

\usepackage{amsmath}
\usepackage{amsthm}
\usepackage{amsfonts}
\usepackage{amssymb}
\theoremstyle{break}
\newtheorem{Thm}{Theorem}
\newtheorem{Lem}{Lemma}
\newtheorem{Cor}{Corollary}
\newtheorem{Rem}{Remark}

\newcommand{\Com}{\mathbb{C}}

\newcommand{\ov}{\overline}
\newcommand{\Real}{\mathbb{R}}

\normalsize

\begin{document}

\title{Analysis of the horizontal Laplacian for the Hopf fibration}
\author{Robert O. Bauer}
\address{Department of Mathematics\\ University of Illinois at Urbana-Champaign\\
 Urbana, IL 61801}
\email{rbauer@math.uiuc.edu}

\keywords{Horizontal Laplacian, Hopf fibration, heat trace asymptotics, Green function, Sobolev inequality}

\subjclass{35K05, 35P10, 58J35}

\begin{abstract}
	We study the horizontal Laplacian $\Delta^H$ associated to 	the Hopf fibration $S^3\to S^2$ with arbitrary Chern  number 	$k$. We use representation theory to calculate the spectrum, 	describe the heat kernel and 	obtain the complete heat 	trace asymptotics of $\Delta^H$.  	We express the Green 	functions for associated Poisson 	semigroups and obtain 	bounds for their contraction properties 	and Sobolev inequalities for $\Delta^H$. The bounds and 
	inequalities improve 	as $|k|$ increases.     
\end{abstract}

\maketitle

\section{Introduction}

In this paper we are concerned with the analysis of the horizontal (or Bochner) Laplacian $\Delta^H$ for the Hopf fibration $S^3\to S^2$. By analysis we mean spectral analysis, description of the heat kernel and its trace, $L^p$ bounds for Green functions of the  associated Poisson semigroup and Sobolev inequalities. On a general level such results are readily available in the literature: For a section $f$ one applies the diamagnetic or Kato's inequality
\[
		\langle \Delta^H f,f\rangle\ge\langle\Delta|f|,|f|\rangle
\]
and now well known estimates and Sobolev inequalities for the scalar Laplacian $\Delta$ translate into estimates and inequalities for the vector Laplacian $\Delta^H$. Of course, the diamagnetic inequality rarely leads to sharp inequalities. The point of this paper  is to improve these estimates and inequalities by directly studying the spectrum of a Laplacian on sections in a particularly simple, yet geometrically non-trivial case,  the Hopf fibration. As it turns out, the curvature of the bundle improves the analytical properties of $\Delta^H$ (vis-a-vis $\Delta$), as indicated for example, by a nonzero first eigenvalue.

The idea for this paper grew out of \cite{bauer.carlen:2001} where the radial part of the spectrum of the horizontal Laplacian on more general Hopf fiberings was determined by studying expectations of the stochastic parallel transport in the bundle. This in turn was motivated by the numerous stochastic approaches to spectral analysis of Schr\"odinger operators with  magnetic fields, for example \cite{ueki:1999}, \cite{erdos:1994}, \cite{ikeda.et.al:1995}, \cite{malliavin:1982} and references therein. Geometrically, the magnetic potential $A$ is a connection 1-form for a trivial complex line bundle over $R^n$, the magnetic field $B=dA$ its curvature, and the corresponding Schr\"odinger operator with magnetic field $-(\nabla+iA)^2$ is the horizontal Laplacian. Often one considers   $\lambda A$ and the spectral properties of the associated Schr\"odinger operator as $\lambda\to\infty$. This corresponds to ``turning up'' the field strength.   In \cite{ueki:1999}, \cite{erdos:1994}, \cite{ikeda.et.al:1995}, \cite{malliavin:1982} the analysis proceeds through the study of certain stochastic oscillatory integrals. In the case considered here we use  representation theory. Note that because we work over $S^2$ there are topological restrictions on our bundle and its connection. This restriction is the Chern number of the bundle. We consider the complete family of Hopf fiberings $S^3\to S^2$ labeled by the Chern number $k\in Z$. Large $|k|$ corresponds to large curvature and thus to a strong ``magnetic field.'' 

We begin by describing the Hopf fibration $S^3\to S^2$ and the associated horizontal Laplacian $\Delta^H$. We calculate its eigenvalues in Lemma~\ref{L:eigenvalues} and its eigenspaces in Lemma~\ref{L:spectrum}. From this we obtain an expression for the heat kernel $e^{-t\Delta^H}$ in Lemma~\ref{L:heatkernel}. Using a formula due to Ramanujan we calculate the heat trace asymptotics in Theorem~\ref{T:asymptotics}. In particular, we obtain a common formula for the asymptotic expansion of $\text{Tr}(e^{-t\Delta^H})$ for Hopf fiberings of all charges $k$. The case $k=0$ corresponds to a trivial bundle where $\Delta^H$ reduces to $\Delta_{S^2}$ and we recover the well known formula  of McKean and Singer. The common formula sheds some light on a question raised by Gilkey: What is the relation between the heat trace asymptotics of Riemannian submersions with the same base and same fiber.

Next, we consider the Green functions for the Poisson semigroup $e^{-t(\Delta^H)^{1/2}}$ and also in the case where an ``appropriate'' mass-term has been added, $e^{-t(\Delta+k^2+1)^{1/2}}$. In Theorem~\ref{T:exgreen} we obtain a closed formula for the Green function with the mass-term. We use this formula to calculate the $L^1$-norm of the Green function, which leads to a bound for the $L^p\to L^p$ norm of the semigroup, stated in Theorem~\ref{T:ck}. The contraction becomes stronger as the charge $|k|$ increases. 

In Theorem~\ref{T:sobolev} and its corollaries we obtain Sobolev inequalities for $\Delta^H$. Here, too, the constants in the inequalities improve with $|k|$. In particular, we show in Corollary~\ref{C:log} that if $\langle\Delta^H f,f\rangle=1$, then $$\ln\int e^{|f|}\ d\mu=o(1),\  \text{as } |k|\to\infty.$$ By comparison, Onofri's inequality---a sharp form of a Moser-Trudinger inequality---combined with the diamagnetic inequality only gives
$$\ln\int_{S^2}e^{|f|}\ d\mu\le\frac{1}{4}+\int_{S^2}|f|\ d\mu.$$ The sharp form of Corollary~\ref{C:log} has eluded us and it would be interesting to see if techniques used to prove Onofri's inequality, such as spherically symmetric rearrangements, see \cite{carlen.loss:1992}, or regularized determinants, see \cite{sarnak.et.al:1988}, can be adapted to obtain a sharp inequality for $\Delta^H$.

\section{Spectrum of the horizontal Laplacian}

\subsection{Laplacian on $S^3$}

Let
$S^3=\{(x^1,x^2,x^3,x^4)\in\Real^4:(x^1)^2+(x^2)^2+(x^3)^2+(x^4)^2=1\}$.
For $1\le\alpha,\beta\le4$ set
\[
    D_{\alpha\beta}=-x^{\alpha}\frac{\partial}{\partial
    x^{\beta}}+x^{\beta}\frac{\partial}{\partial x^{\alpha}},
\]
and let
\[
    K_x=D_{14}+D_{23},\quad K_y=D_{31}+D_{24},
    \quad K_z=D_{12}+D_{34}.
\]
The Laplace-Beltrami operator on $S^3$ with its standard metric is
given by
\[
    \Delta=-(K_x^2+K_y^2+K_z^2),
\]
see \cite[Section 1.2]{englefield:72}.


\subsection{Hopf fibration $S^3\to S^2$ }

Let $a=x^1+ix^2$, $b=x^3+ix^4$. Then $S^3=\{(a,b)\in\Com^2:
a\ov{a}+b\ov{b}=1\}$ and $S^1=\{u\in\Com:u\ov{u}=1\}$ acts on
$S^3$ by
\[
    (a,b)u=(au,bu).
\]
This action induces an $S^1$-principal bundle over the quotient
manifold $P_1\Com\simeq S^2$, which is known as the Hopf fibration
\[
    S^1\to S^3\to S^2.
\]
For the natural connection on this principal bundle the vector
fields $K_x$ and $K_y$ are horizontal. The vector field $-K_z$
induces the $S^1$-action and is vertical. The horizontal Laplacian
for the Hopf fibration is
\[
    \Delta^H=-(K_x^2+K_y^2).
\]
Given a representation $\rho$ of $S^1$ on $\Com$ we obtain a
complex line bundle $E\to S^2$ associated to the $S^1$-principal
bundle by setting
\[
    E=S^3\times_{S^1}\Com,
\]
where the $S^1$-action on $S^3\times\Com$ is given by
\[
    ((a,b),z)u=((a,b)u,\rho(u)^{-1}z),\quad (a,b)\in S^3,z\in\Com.
\]
Sections $f$ of the line bundle $E\to S^2$ are given by functions
$f:S^3\to\Com$ which are $\rho$-equivariant
\[
    f((a,b)u)=\rho(u)^{-1}f((a,b)).
\]
Denote $\Gamma(E)$the space of sections of $E\to S^2$. We are
interested in the operator $\Delta^H$ acting on $\Gamma(E)$.


\subsection{Spectrum of $\Delta^H|_{\Gamma(E)}$ } Identify $S^1$
with the toral group $\{e^{it}:t\in\Real\}$. For integer $k$
define a 1-dimensional irreducible representation $\rho_k$ by
\[
    \rho_k(e^{it})=e^{ikt}.
\]
Assume now that $f\in\Gamma(E_{k})$, where $E_k$ is the complex
line bundle obtained from the representation $\rho_k$. Then $f$ is
an eigenfunction of $K_z$ with eigenvalue $ik$. Indeed,
\begin{align}
    K_z
    f(a,b)&=\lim_{s\to0}\frac{f((a,b)e^{-is})-f((a,b))}{s}\notag\\
    &=f((a,b))\lim_{s\to0}\frac{e^{iks}-1}{s}=ik f((a,b)).
\end{align}
Thus, if $f\in\Gamma(E_k)$ is an eigenfunction of $\Delta^H$ with
eigenvalue $\lambda$, then
\[
    \Delta f=(\Delta^H-K_z^2)f=(\lambda+k^2)f.
\]
The eigenvalues of $\Delta$ are $N^2-1$, $N=1,2,\dots$. Set
$N=2l+1$. The eigenspace for the eigenvalue $N^2-1$ has dimension
$N^2$ and contains an $N$-dimensional subspace of simultaneous
eigenfunctions of $K_z$ with eigenvalue $2ip$. Here $p$ is such
that $l+p$ is integer and $-l\le p\le l$, see \cite[Section
1.3]{englefield:72}. Thus $k/2$ and $l$ are either both integer or
both non-integer. Hence $N-1$ and $k$ have equal parity.
Furthermore, $l\ge |k/2|$ and so $N\ge|k|+1$. Thus $N=|k|+1+2M$,
where $M=0,1,2,\dots$ and so

\begin{Lem} \label{L:eigenvalues}
The eigenvalues of $\Delta^H$ on
$\Gamma(E_k)$ are
\begin{equation}\label{E:eigenvalue}
    N^2-1-k^2=4M^2+4(1+|k|)M+2|k|, \quad M=0,1,2,\dots.
\end{equation}
In particular, $\Delta^H$ has a mass gap of size $2|k|$.
\end{Lem}


\subsection{Eigenfunctions for $\Delta^H|_{\Gamma(E_k)}$}

Identify $S^3$ with $SU(2)$ by
\[
    (a,b)\mapsto\begin{bmatrix}a&
    b\\-\ov{b}&\ov{a}\end{bmatrix},\quad a\ov{a}+b\ov{b}=1.
\]
For $g\in SU(2)$ set
\begin{equation}\label{E:homogeneous}
    t_{mn}^l(g)=\sqrt{\frac{(l+m)!}{(l-n)!(l+n)!(l-m)!}}
    \frac{d^{l-m}}{dx^{l-m}}[(ax-\ov{b})^{l-n}(bx+\ov{a})^{l+n}]|_{x=0},
\end{equation}
where $l$ is an integer multiple of 1/2, $m$ and $n$ are such that
$l+m$ and $l+n$ are integer, and $-l\le n,m\le l$. As shown in
\cite[Section 9.14]{askey:1999}, the $t_{mn}^l$ are matrix entries
of an irreducible unitary $N=2l+1$-dimensional representation. A
simpler form for the $t_{mn}^l$ can be given in terms of Euler
angles. The Euler angles $\phi$, $\psi$, and $\theta$ are obtained
from the relations
\[
    a=e^{i(\phi-\psi)/2}\cos{\frac{1}{2}\theta},\quad
    b=ie^{i(\phi-\psi)/2}\sin{\frac{1}{2}\theta},
\]
where $0\le\phi<2\pi,\ 0<\theta<\pi,$ and $-2\pi\le\psi<2\pi$.
Note that
\begin{equation}\label{E:matrix}
    \begin{bmatrix}a& b\\-\ov{b}&\ov{a}\end{bmatrix}=
    \begin{bmatrix}e^{i\phi/2}&0\\0& e^{-i\phi/2}\end{bmatrix}
    \begin{bmatrix}\cos{\frac{1}{2}\theta}&
    i\sin{\frac{1}{2}\theta}\\
    i\sin{\frac{1}{2}\theta}&\cos{\frac{1}{2}\theta}\end{bmatrix}
    \begin{bmatrix}e^{i\psi/2}&0\\0&e^{-i\psi/2}\end{bmatrix}.
\end{equation}
Then
\begin{equation}\label{E:tlmn}
    t_{mn}^l(g)=e^{-i(m\phi+n\psi)}P_{mn}^l(\cos\theta),
\end{equation}
where
\begin{align}
    P_{mn}^l(z)=&\frac{(-1)^{l-n}i^{n-m}}{2^l}
    \sqrt{\frac{(l+m)!}{(l-n)!(l+n)!(l-m)!}}(1+z)^{-(m+n)/2}\notag\\
    &\cdot(1-z)^{(n-m)/2}\frac{d^{l-m}}{dz^{l-m}}
    [(1-z)^{l-n}(1+z)^{l+n}].
\end{align}
Note that $P_{mn}^l(z)$ can be written in terms of a Jacobi
polynomial. It s a constant multiple of
\[
    (1-z)^{(m-n)/2}(1+z)^{(m+n)/2}P_{l-m}^{(m-n,m+n)}(z).
\]
From \eqref{E:homogeneous} it is easy to see that $t_{mn}^l$ is an
eigenfunction of $\Delta$ with eigenvalue $N^2-1$. Furthermore,
from \eqref{E:matrix} wee see that if
\[
    (\phi,\psi,\theta)\equiv(a,b),\quad
    (\phi',\psi',\theta')\equiv(au,bu),\ u=e^{it},
\]
then
\[
    \phi'=\phi+2t,\ \psi'=\psi,\ \theta'=\theta.
\]
Using $K_z=-2\partial/\partial\phi$ and \eqref{E:tlmn}, we get
\[
    K_z t_{mn}^l(g)=2im\ t_{mn}^l(g).
\]
Thus

\begin{Lem}\label{L:spectrum}
	For $M=0,1,2\dots$, the eigenspace of 	$\Delta^H|_{\Gamma(E_k)}$ for the eigenvalue 	$4M^2+4(1+|k|)M+2|k|$ is spanned by the 
	$2l+1=2M+1+|k|$ eigenfunctions
	$t_{\frac{k}{2},n}^l$, $-l\le n\le l$.
\end{Lem}

\begin{Rem} 
	The functions $t_{mn}^l$ are not normalized:
	\[
		\int_{SU(2)}|t_{mn}^l(g)|^2\ dg=\frac{1}{2l+1},
	\] 
	where $dg$ denotes Haar measure on $SU(2)$, see 	\eqref{E:haar}.
\end{Rem}


\section{Heat kernel, Green function, and Sobolev inequality}

\subsection{Heat kernel, addition formula, and asymptotics}

The heat semigroup $e^{-t\Delta^H}$ on $\Gamma(E_k)$ has a kernel
$k_t(g,g')$ satisfying
\[
    e^{-t\Delta^H}f(g)=\int_{S^3}k_t(g,g')f(g')\ dg'.
\]
The kernel is given explicitly by
\begin{equation}\label{E:explicitkernel}
    	k_t(g,g')=\sum_{M=0}^{\infty} 
	e^{-t\lambda_M}(2l+1)\sum_{n=-l}^l
    	t_{\frac{k}{2},n}^l(g)\ \ov{t_{\frac{k}{2},n}^l(g')},
\end{equation}
where $\lambda_M$ is given by \eqref{E:eigenvalue}, see
\cite[Section 1.6]{gilkey:1999}. 

\begin{Lem}\label{L:heatkernel}
	We have
	\begin{equation}\label{E:heatkernel}
    		k_t(g,g')=\sum_{M=0}^{\infty}e^			{-t(4M^2+4(1+|k|)M+2|k|)}(|k|+2M+1)
    		t_{\frac{k}{2},\frac{k}{2}}^{M+|k|/2}(g(g')^{-1}).
\end{equation}
\end{Lem}

\begin{proof}
	Since the $t_{mn}^l(g)$ are matrix
	entries of a unitary representation we have
	\begin{align}\label{E:addition}
    		\sum_{n=-l}^l t_{mn}^l(g)\ \ov{t_{mn}^l(g')}&=\sum_{n=-l}^l
    		t_{mn}^l(g)\ t_{nm}^l((g')^{-1})\notag\\ &=t_{mm}^l(g(g')^		{-1}).
	\end{align}
	The lemma follows.
\end{proof}

From this lemma we get 
\begin{equation}\label{E:kerneltrace}
	\text{Tr}_{L^2}(e^{-t\Delta^H})=\sum_{M=0}^{\infty}e^			{-t(4M^2+4(1+|k|)M+2|k|)}(|k|+2M+1),
\end{equation}
see  \cite{mckean.singer:1967}, and \cite[Section 1.6.12]{gilkey:1999}. Using a formula due to Ramanujan we can calculate the small $t$ behavior of this trace.

\begin{Thm}[Heat equation asymptotics]\label{T:asymptotics}
	The trace of the heat kernel $e^{-t\Delta^H}$ on sections of 	$\Gamma(E_k)$ has the asymptotic expansion
	\[
		2e^{t(k^2+1)}\left(\frac{1}{8t}+\sum_{r=0}^{\infty} 		\zeta((1+|k|)/2,-1-2r)\frac{(-4t)^r}{r!}\right),
	\]	
	where the Hurwitz zeta function is given by
	\[
		\zeta(x,s)=\sum_{n=0}^{\infty}\frac{1}{(n+x)^s}.
	\]
	In particular  
	\begin{align}\label{E:expand5}
		\text{Tr}_{L^2}(e^{-t\Delta^H})=&\frac{1}{4t}
		+\frac{1}{3}\notag\\
		&+\frac{1}{30}(8-5k^2)t \notag\\
		&+\frac{1}{315}(64-126k^2)t^2\notag \\ 
		&+\frac{1}{630}(128-432k^2+49k^4)t^3\notag\\
		&+\frac{2}{3465}(512-2112k^2+561k^4)t^4\notag\\
		&+\frac{1}{675675} 							(391168-1722240k^2+669240k^4-22165k^6)t^5
		+O(t^6).
	\end{align}
\end{Thm}

\begin{proof}
	From \eqref{E:kerneltrace} we have
	\[
		\text{Tr}_{L^2}(e^{-t\Delta^H}) =2e^{t(k^2+1)}		\sum_{M=0}^{\infty}e^{-4t(M+(1+|k|)/2)^2}
		(M+(1+|k|)/2).
	\]
	Now apply Lemma  \ref{L:asymptotic} from the appendix. 
	A less 	concise 	expression for the expansion uses 
	\begin{align}\label{E:expandcalc}
		&\sum_{M=0}^{\infty}e^{-4t(M+(1+|k|)/2)^2}
		(M+(1+|k|)/2)\notag\\
		&=\frac{1}{8t}+\sum_{r=0}^{\infty}\zeta(-1-2r) 
		\frac{(-4t)^r}{r!} 		
		-\sum_{M=1}^{(|k|-1)/2}e^{-4tM^2}M,\ \text{for odd }k 		\notag\\
		&=\frac{1}{8t}+\sum_{r=0}^{\infty}\zeta(\frac{1}{2},-1-2r) 		\frac{(-4t)^r}{r!}-\sum_{M=0}^{|k|/2-1}e^{-4t(M+1/2)^2} 		(M+\frac{1}{2}),\ \text{for even }k.
	\end{align}
	Here $\zeta(s)=\zeta(1,s)$ is the Riemann zeta function. This 	expression is derived in the same way as Lemma 	\ref{L:asymptotic} in the appendix.
	The coefficients up to order five are calculated from this 	expression using Mathematica.
\end{proof}

\begin{Rem}
	Note that for any $k$ the first two terms are  	$1/(4t)+1/3$, expressing the dimension and topology 	of the base. The bundle is only visible in the higher order 	terms.
	Note also that for $k=0$ the horizontal Laplacian reduces to 	the 	Laplacian on functions on $S^2$ with the round sphere 	metric 	of radius $1/2$ and the volume normalized to 1. We 	recover 	the formula from \cite[Table,  	p.63]{mckean.singer:1967}
	\[
		\text{Tr}_{L^2}(e^{-t\Delta})=\frac{1}{4t}+\frac{1}{3} 		+\frac{4t}{15}+\frac{64t^2}{315}+\frac{64t^3}{315} 		+\frac{1024t^4}{3465}+\frac{391168t^5}{675675}+O(t^6). 
	\]
	Because McKean and Singer work with the sphere of radius 1 	and volume $4\pi$ one has to substitute 1 for $4\pi$ and $4t$ 	for $t$ to obtain our formula from theirs.
\end{Rem}

\begin{Rem}
	By the preceding remark Theorem \ref{T:asymptotics} 	relates the heat trace 	asymptotics of the horizontal 	Laplacian for the standard Hopf 	fibration (i.e. $k=1$) and of 	the skalar Laplacian on $S^2$ 	through a common formula. 	The 	problem of relating heat trace 	asymptotics of Riemannian 	submersions with same base 	space was 	raised in \cite[Section 	4.7]{gilkey:1999}.  	According to  \eqref{E:expand5}, explaining the relation 
	means explaining what the coefficients of $k^2, 	k^4,k^6,\dots$ in the respective coefficients of $t$ in 	\eqref{E:expand5} stand for.
\end{Rem}       


\subsection{Green functions}

Denote $h_t(g,g')$ and $h_t^{\#}(g,g')$ the kernels of the Poisson
semigroups $e^{-t(\Delta^H)^{1/2}}$ and $
e^{-t(\Delta^H+k^2+1)^{1/2}}$, respectively. Since all eigenvalues
of $\Delta^H+k^2+1$ are squares of integers, the formulas for
$h_t^{\#}$ will be simpler. Note that $k^2+1$ can be considered the square of the mass of the particle described by the operator $\Delta^H+k^2+1$. From \eqref{E:heatkernel},
\begin{equation}\label{E:poissonkernel}
    h_t(g,g')=\sum_{M=0}^{\infty}e^{-t((|k|+2M+1)^2-1-k^2)^{1/2}}(|k|+2M+1)
    t_{\frac{k}{2},\frac{k}{2}}^{M+|k|/2}(g(g')^{-1}),
\end{equation}
and
\begin{equation}\label{E:poissonkernelsharp}
    h_t^{\#}(g,g')=\sum_{M=0}^{\infty}e^{-t(|k|+2M+1)}(|k|+2M+1)
    t_{\frac{k}{2},\frac{k}{2}}^{M+|k|/2}(g(g')^{-1}).
\end{equation}
Integration over $t\in(0,\infty)$ gives the Green functions
\begin{equation}\label{E:green1}
    G(g,g')=\sum_{M=0}^{\infty}
    \frac{|k|+2M+1}{((|k|+2M+1)^2-1-k^2)^{1/2}}
    t_{\frac{k}{2},\frac{k}{2}}^{M+|k|/2}(g(g')^{-1}),
\end{equation}
and
\begin{equation}\label{E:green2}
    G^{\#}(g,g')=\sum_{M=0}^{\infty}
    t_{\frac{k}{2},\frac{k}{2}}^{M+|k|/2}(g(g')^{-1}).
\end{equation}

\begin{Thm}\label{T:exgreen}
The Green function
$G^{\#}$ is given by
\begin{equation}\label{E:green}
    G^{\#}(g,g')=e^{-i|k|(\phi+\psi)/2}\frac{1}{2\sin\frac{\theta}{2}}
    \left(\frac{\cos\frac{\theta}{2}}{1+\sin\frac{\theta}{2}}\right)^{|k|},
\end{equation}
where $g(g')^{-1}=(\phi,\psi,\theta)$ if $k>0$ and \begin{equation}\label{E:minusk}
    \begin{bmatrix}0& i\\i&0\end{bmatrix}\circ g(g')^{-1}\circ 	\begin{bmatrix}0& i\\i&0\end{bmatrix}=(\phi,\psi,\theta)
\end{equation}
if $k<0$.
\end{Thm}

\begin{proof}
As mentioned above, we can express the eigenfunctions in terms of Jacobi polynomials,
\[
    t_{\frac{k}{2},\frac{k}{2}}^{M+|k|/2}(g(g')^{-1})
    =e^{-ik(\phi+\psi)/2}2^{-k/2}(1+z)^{k/2}
    P_{M+(|k|-k)/2}^{(0,k)}(z),
\]
 where $z=\cos\theta$. The
generating function for Jacobi polynomials
\[
    2^{\alpha+\beta}R^{-1}(1-r+R)^{-\alpha}(1+r+R)^{-\beta}
    =\sum_{n=0}^{\infty}P_n^{(\alpha,\beta)}(x)r^n,
\]
where $R=(1-2xr+r^2)^{1/2}$, together with \eqref{E:green2} then imply the result for $k>0$. Now note that
\[
	\left[t^l\left(\begin{bmatrix}0& i\\i&0\end{bmatrix}\right)\ t^l(g)\ 	t^l\left(\begin{bmatrix}0& i\\i&0\end{bmatrix}\right) 	\right]_{mn}=\left[t^l(g)\right]_{-m,-n}.
\]
Thus results for the $(k,k)$-entry of the matrix $t^l$ translate into results for the $(-k,-k)$-entry. This proves the theorem for $k<0$. 
\end{proof} 

\subsection{Inequalities}

The Haar measure on $SU(2)$ is given by
\begin{equation}\label{E:haar}
    dg=\frac{1}{16\pi^2}\sin\theta\ d\theta d\phi d\psi.
\end{equation}
Set
\begin{equation}\label{E:g1}
    c_k=\int_{SU(2)}|G^{\#}(g,g')|dg'.
\end{equation}
Note that because $G^{\#}(g,g')$ only depends on $g(g')^{-1}$, $c_k$ does not depend on $g$. 

\begin{Lem}
	If $k$ is even, $|k|=2m$, $m=1,2,\dots$, then
	\[
		c_k=2m|\ln 2-\sum_{n=1}^{m-1}(-1)^{n-1}/n|-1
	\]
	and if $k$ is odd, $|k|=2m+1$, $m=0,1,\dots$, then
	\[
		c_k=2(2m+1)|\frac{\pi}{4}-\sum_{n=0}^{m-1}		(-1)^n/(2n+1)|-1,
	\]
	where the empty sums are taken to be zero. Furthermore,
	\begin{equation}\label{E:ckineq}
		\frac{1}{|k|+2}< c_k <\frac{1}{|k|},\quad |k|=1,2,\dots
	\end{equation}
\end{Lem}

\begin{proof} We may assume that $k>0$.
	From \eqref{E:green} follows
	\[
		c_k=\frac{1}{2}\int_0^{\pi}
    		\left(\frac{\cos(\theta/2)}{1+\sin(\theta/2)}\right)^k
		\cos(\theta/2)\ d\theta.
	\]
	Substituting $\theta$ for $\theta/2$, using 	$\sin\theta=\cos(\pi/2-\theta)$ and 	$\tan(\theta/2)=\sin\theta/(1+\cos\theta)$ gives
	\[
		c_k=4\int_0^{\pi/4}\tan^k\theta\sin\theta\cos\theta\ 		d\theta.
	\]
	Finally, the substitution $x=\tan\theta$ gives 
	\begin{equation}\label{E:ckint}
		c_k=4\int_0^1\frac{x^{1+k}}{(1+x^2)^2}\ dx.
	\end{equation}
	Note that
	\[
		x^{1+k}\le\frac{4x^{1+k}}{(1+x^2)^2}\le x^{k-1}
	\]
	for $0\le x\le1$. Since the inequalities are strict for $x\neq 0,1$ 	\eqref{E:ckineq} follows.  Integrating by parts in \eqref{E:ckint} 	and power series gives 
	\[
		c_k=2k\sum_{n=0}^{\infty}\frac{(-1)^n}{2n+k}\ -1.
	\]
	The explicit formulae now follow from 
	\[
		\sum_{n=1}^{\infty}\frac						{(-1)^{n-1}}{n}=\ln2,\quad\sum_{n=0}^{\infty}
		\frac{(-1)^{n}}{2n+1}=\frac{\pi}{4}.
	\]
\end{proof}
In particular
\[
	c_1=\frac{\pi}{2}-1,\quad c_2=2\ln 2-1.
\]

\begin{Thm}\label{T:ck}
	For a section $f\in\Gamma(E_k)$ we have
	\begin{equation}\label{E:gsharp}
		 \|G^{\#}f\|_p\le c_k\|f\|_p.
	\end{equation}
\end{Thm}

\begin{proof}
By H\"older's inequality we have
\[
    |\int G^{\#}(g,g')f(g')\ dg'|^p\le(\int|G^{\#}(g,g')|\
    dg')^{p-1}\int|G^{\#}(g,g')||f(g')|^p\ dg'.
\]
The result follows by applying Fubini-Tonelli and \eqref{E:g1}. 
\end{proof}

\begin{Cor}
	The inequality \eqref{E:ckineq} can be improved to
	\begin{equation}\label{E:ckineqplus}
		\frac{1}{1+|k|}\le c_k<\frac{1}{|k|}.
	\end{equation}
\end{Cor}

\begin{proof}
	Note that the first eigenvalue of $(\Delta^H+1+k^2)^{1/2}$ is 	$1+|k|$. Thus for $p=2$ we get from the 	eigenfunction representation of $G^{\#}$ immediately the 	sharp bound
	\[
    		\|G^{\#}f\|_2\le \frac{1}{1+|k|}\|f\|_2.
	\]
	Indeed, for a section $f$ write
	\[
    		f(g)=\sum_{M=0}^{\infty}\sqrt{2l+1}\sum_{n=-l}^l  a_{ln}
    		t_{\frac{k}{2},n}^l(g),
	\]
	and note that directly from the expression 	\eqref{E:explicitkernel} 
	\[
		G^{\#}(g,g')=\sum_{M=0}^{\infty} \sum_{n=-l}^l
    		t_{\frac{k}{2},n}^l(g)\ \ov{t_{\frac{k}{2},n}^l(g')}.
	\]
	Thus
	\[
		G^{\#}f(g)=\sum_{M=0}^{\infty} 			\sum_{n=-l}^l\frac{1}{\sqrt{2l+1}} a_{ln} t^l_{\frac{k}{2},n}(g)
	\]
	and
	\[
		\|G^{\#}f\|_2^2=\sum_{M=0}^{\infty} 			\sum_{n=-l}^l\frac{1}{(2l+1)^2}|a_{ln}|^2.
	\]
	For given $\|f\|_2^2=\sum_{M=0}^{\infty} \sum_{n=-l}^l 	|a_{ln}|^2$ this expression is maximal if $a_{ln}=0$ for $M>0$. 	In that case 
	\[
		\|G^{\#}f\|_2=\frac{1}{1+|k|}\|f\|_2.
	\]
	The Corollary now follows from Theorem \ref{T:ck}.
\end{proof}

\begin{Rem}
	For $p=2$ we can also calculate $\|Gf\|_2$. Indeed, from 	\eqref{E:explicitkernel} and \eqref{E:green1}
	\[ 
		G(g,g')=\sum_{M=0}^{\infty} 				\frac{|k|+2M+1}{((|k|+2M+1)^2-1-k^2)^{1/2}} \sum_{n=-l}^l
    			t_{\frac{k}{2},n}^l(g)\ \ov{t_{\frac{k}{2},n}^l(g')}.
	\]
	Thus
	\[
		Gf(g)=\sum_{M=0}^{\infty} 					\left(\frac{|k|+2M+1}{(|k|+2M+1)^2-1-k^2}\right)^{1/2} 		\sum_{n=-l}^l  a_{ln} t^l_{\frac{k}{2},n}(g)
	\]
	and
	\[
		\|Gf\|_2^2=\sum_{M=0}^{\infty} 			\frac{1}{(|k|+2M+1)^2-1-k^2} 
		\sum_{n=-l}^l |a_{ln}|^2.
	\]
	As in the proof above, this expression is maximal if $a_{ln}$ is 	zero for $M>0$. Thus
	\begin{equation}\label{E:greenL2}
		\|Gf\|_2^2\le \frac{1}{2|k|}\|f\|_2^2
	\end{equation}
	and
	\begin{equation}\label{E:greenL1}
		\|Gf\|_1\le\frac{1}{\sqrt{2|k|}}\|f\|_2.
	\end{equation}

\end{Rem}
	
\begin{Thm}\label{T:sobolev}[Sobolev Inequality]
	For a section $f\in\Gamma(E_k)$ with $\|f\|_2=1$ and 	$p>2$ we have
	\begin{equation}\label{E:sob1}
		\|Gf\|_p^p\le|k|^								{-1/3}\left(1+2^{(1-p)/(p-2)}(p-2)\right)^{(p-2)/2},
	\end{equation}
	and
	\begin{equation}\label{E:sob2}
		\|G^{\#}f\|_p^p\le|k|^{-1}\left(\frac{p-2}{2}\right)^{(p-2)/2}.
	\end{equation}
\end{Thm}

\begin{proof}
	Since $\|f\|_2=1$, we have $\sum_{l,n} |a_{ln}|^2=1$. Set
	$a_l=(\sum_{n=-l}^l |a_{ln}|^2)^{1/2}$. Then 	$\sum_{l}a_l^2=1$, and 	$\frac{1}{a_l}(a_{l,-l},a_{l,-l+1},\dots,a_{l,l})$ is a unit
	vector. Thus there exists a $g_l^{'}\in SU(2)$ so that
	\[
    		\frac{a_{ln}}{a_l}=t_{\frac{k}{2},n}^l(g_l^{'}),\quad -l\le n\le
    		l.
	\]
	Since the $t_{mn}^l$ are matrix entries of a unitary
	representation it follows from \eqref{E:greeneigen} that
	\[
    		Gf(g)=\sum_{M=0}^{\infty}\frac{\sqrt{2l+1}}{\lambda_M}a_l
    		t_{\frac{k}{2},\frac{k}{2}}^l(g_l),
	\]
	for some $g_l\in SU(2)$. Writing $a_M$ for $a_l$, 	$t_M(g_M)$ for $t_{\frac{k}{2},\frac{k}{2}}^l(g_l)$, and setting
	\begin{equation}\label{E:bm}
		b_M=\frac{\lambda_M^2}{2l+1}=2M+|k|+1- 		\frac{k^2+1}{2M+|k|+1} 
	\end{equation}
	we get the expression
	\[
		Gf(g)=\sum_{M=0}^{\infty}b_M^{-1/2}a_M t_M(g_M).
	\]
	From H\"older's inequality
	\[
		|Gf(g)|^p\le\left(\sum_{M=0}^{\infty}b_M^{-1/2} 		a_M^{p/(p-1)}\right)^{p-1}\sum_{M=0}^{\infty}b_M^{-1/2} 		|t_M(g_M)|^p.
	\]
	Applying H\"older again we get
	\[
		|Gf(g)|^p\le\left(\sum_{M=0}^{\infty} 			b_M^{(1-p)/(2-p)}\right)^{(p-2)/2} 			\left(\sum_{M=0}^{\infty}a_M^2\right)^{p/2} 		\sum_{M=0}^{\infty}b_M^{-1/2}|t_M(g_M)|^p.
	\]
	The middle sum equals 1 by assumption. As entries of a 	unitary matrix $|t_M(g_M)|\le1$ and so 
	\[
		\int_{SU(2)}|t_M(g_M)|^p\ dg_M\le\int_{SU(2)}|t_M(g_M)|^2\ 		dg_M=\frac{1}{2M+|k|+1}.
	\]
	Thus, after integration, the third sum is bounded by 	$\sum b_M^{-1/2}(2M+|k|+1)^{-1}$. From \eqref{E:bm} we get 	\[
		b_M\ge2M+1+\frac{|k|-1}{|k|+1}.
	\] 
	Since
	\[
		\left(2M+1+\frac{|k|-1}{|k|+1}\right)\left(2M+1+|k|\right)^2\ge 		\left(2M+1+|k|^{2/3}\right)^3,
	\]
	the third sum is bounded above by
	\[
		\sum_{M=0}^{\infty}\left(2M+1+|k|^{2/3}\right)^{-3/2}.
	\]
	By simple integral comparison we now get
	\begin{equation}\label{E:third}
		\sum_{M=0}^{\infty}b_M^{-1/2}|t_M(g_M)|^p\le|k|^{-1/3}.
	\end{equation}
	Similarly, using $b_M\ge 2M+1$, we get
	\begin{equation}\label{E:first}
		\sum_{M=0}^{\infty}b_M^{(1-p)/(p-2)}
		\le1+2^{(1-p)/(p-2)} (p-2)
	\end{equation}
	Combining \eqref{E:third} and \eqref{E:first} gives 	\eqref{E:sob1}.
	For $G^{\#}$ we have 
	\[
		b_M=2M+1+|k|
	\]
	and so 
	\[
		\sum_{M=0}^{\infty}b_M^{(1-p)/(p-2)}\le|k|^		{-1/(p-2)}\left(\frac{p-2}{2}\right).
	\]
	The third sum above in this case agrees with the first sum  	with $p=4$. Combining this gives \eqref{E:sob2}.
\end{proof}

\begin{Cor}
	For a section $f\in\Gamma(E_k)$ and $p>2$ we have
	\[
		\|f\|_p\le |k|^								{-1/3p}\left(1+2^{(1-p)/(p-2)}(p-2)\right)^{(p-2)/2p} 	\langle\Delta^H f,f\rangle^{1/2}.
	\]
\end{Cor}

Note that for any $p$ we can achieve $\|f\|_p\le\langle\Delta^H f,f\rangle^{1/2}$ by choosing $|k|$ sufficiently large.

\begin{Cor}\label{C:log}
	For a section $f\in\Gamma(E_k)$ with $\|f\|_2=1$ there are  	constants $C$ and $C^{\#}$ so that
	\[
		\ln\left(\int_{SU(2)}e^{|Gf(g)|}\ dg\right)\le|k|^{-1/3}C,
	\]
	and 
	\[
		\ln\left(\int_{SU(2)}e^{|G^{\#}f(g)|}\ dg\right)\le|k|^{-1}C^{\#}.
	\]
	\end{Cor}

\begin{proof}
	Note that from \eqref{E:greenL2}, \eqref{E:greenL1} 	we have $\|Gf\|_1\le2^{-1/2}|k|^{-1/3}$ and 	$\|Gf\|_2^2\le\frac{1}{2}|k|^{-1/3}$. Using the 	power series for $e^x$ together with \eqref{E:sob1} and  	$\ln(1+x)\le x$ implies the first inequality with, for example,  
	\[
		C=1+\sum_{n=3}^{\infty}\frac{1}{n!} 		\left(1+2^{(1-n)/(n-2)}(n-2)\right)^{(n-2)/2}.
	\]
	The second inequality follows similarly now using in 	addition the bounds 	from \eqref{E:gsharp}. We can take
	\[
		C^{\#}=c_k+\frac{1}{8k^2}+\sum_{n=3}^{\infty}\frac{1}{n!} 		\left(\frac{n-2}{2}\right)^{(n-2)/2}.
	\]
\end{proof}

\begin{Rem}
	The above Corollary should be compared with Onofri's 	inequality, a sharp form of a 
	Moser-Trudinger inequality,
	\[
		\frac{1}{4}\int_{S^2}|\nabla u|^2\ d\mu+ \int_{S^2}u\ d\mu 		\ge\ln\left(\int_{S^2} e^u\ d\mu\right),
	\]
	where $u$ is a positive function on $S^2$ and $\mu$ is the 	normalized uniform surface measure of $S^2$ 	\cite{carlen.loss:1992}. From the diamagnetic inequality 	\cite{lieb.loss:1996} we get
	\[
		\langle \Delta^H f,f\rangle\ge\langle\Delta |f|,|f|\rangle.
	\]
	Thus Onofri's inequality implies that
	\[
		\frac{1}{4}\int_{S^2}\langle\Delta^H f,f\rangle\ d\mu+ 		\int_{S^2}|f|\ d\mu \ge\ln\left(\int_{S^2} e^{|f|}\ d\mu\right),
	\]
	and, replacing $f$ by $Gf$,
	\[
		\frac{1}{4}\int_{S^2}\langle f,f\rangle\ d\mu+ 		\int_{S^2}|Gf|\ d\mu \ge\ln\left(\int_{S^2} e^{|Gf|}\ 		d\mu\right).
	\]
	Most notable is the absence of a term  	corresponding to $\int_{S^2}u\ d\mu$ in Corollary~\ref{C:log}.  	Its appearance in Onofri's inequality is required because 	addition of  a constant does not change the Dirichlet integral 	of $u$. For sections, we cannot add constants.   
\end{Rem}
	
\appendix
\section{An asymptotic expansion}

\begin{Lem}\label{L:asymptotic}
	Define
	\[
		f(x)=\sum_{M=0}^{\infty}e^{-x(M+\frac{1+|k|}{2})^2} 		(M+\frac{1+|k|}{2}).
	\]
	Then, as $x$ tends to $0+$,
	\[
		f(x)\sim\frac{1}{2x}+\sum_{r=0}^{\infty} 		\zeta((1+|k|)/2,-1-2r)\frac{(-x)^r}{r!}.
	\]
\end{Lem}

\begin{proof}
	The proof proceeds as in \cite{berndt:1985}. For a complex 	number $s$ write $s=\sigma+it$. Using the 	definition of $f$ and inverting the order of summation and 	integration by absolute convergence, we find that 
	\[
		\int_0^{\infty} f(x) x^{s-1}\ dx 				=\Gamma(s)\zeta((1+|k|)/2,2s-1),  
	\]
	provided that $\sigma>1$. Note that for $k$ even, $|k|=2m$, 	$m=0,1,2,\dots$ 
	\begin{equation}\label{E:zeta1}
		\zeta((1+|k|)/2,2s-1) =\zeta(1/2,2s-1) 		-\sum_{r=0}^{m-1}(r+1/2)^{1-2s},
	\end{equation}
	and for $k$ odd, $|k|=2m+1$, $m=0,1,2,\dots$
	\begin{equation}\label{E:zeta2}
		\zeta((1+|k|)/2,2s-1)
		=\zeta(2s-1)-\sum_{r=0}^{m-1}(r+1)^{1-2s},
	\end{equation}
	where $\zeta(s)$ is the Riemann zeta function and empty sums 	are taken to be zero. Note also that 
	\begin{equation}\label{E:zeta3}
		\zeta(1/2, s)=(2^s-1)\zeta(s).
	\end{equation}
	By Mellin's inversion formula 	\cite[p. 33]{titchmarsh:1951}
	\begin{equation}\label{E:inversemellin}
		f(x)=\frac{1}{2\pi i}\int_{2-i\infty}^{2+i\infty} 		\Gamma(s)\zeta((1+|k|)/2,2s-1) x^{-s}\ ds.
	\end{equation}
	Consider now
	\[
		I_{N,T}=\frac{1}{2\pi i}\int_{C_{N,T}} 		\Gamma(s)\zeta((1+|k|)/2,2s-1) x^{-s}\ ds,
	\]
	where $C_{N,T}$ is the positively oriented rectangle with 	vertices $2\pm iT$ and $-(N+\frac{1}{2})\pm iT$, where $T>0$ 	and 	$N>1$ is integer. The integrand has simple poles at 	$s=1$ and  	$s=0,-1,-2,\dots,-N$ in the interior of $C_{N,T}$. 	Note that 	\eqref{E:zeta1}--\eqref{E:zeta3} imply that  	$\zeta((1+|k|)/2,2s-1)$ has the same residue at $s=1$ as 	$\zeta(2s-1)$ and so by the residue theorem
	\begin{equation}\label{E:residue}
		I_{N,T}=\frac{1}{2x}+\sum_{r=0}^N\zeta((1+|k|)/2,-2r-1) 		\frac{(-x)^r}{r!}.
	\end{equation}
	Thus, in order to establish the lemma, it suffices to show that 
	\begin{equation}\label{E:bound1}
		\int_{-(N+\frac{1}{2})}^2\Gamma(\sigma\pm iT) 		\zeta((1+|k|)/2,2(\sigma\pm iT)-1)x^{-\sigma\mp iT}\ 		d\sigma = o(1)
	\end{equation}
	as $T$ tends to $\infty$, and then that
	\begin{equation}\label{E:bound2}
		\int_{-\infty}^{\infty}\Gamma(-(N+\frac{1}{2})+it) 		\zeta((1+|k|)/2,2(-(N+\frac{1}{2})+it)-1) x^{N+\frac{1}{2}-it}\ 		dt\ll x^{N+\frac{1}{2}},
	\end{equation}
	as $x\to0+$. 

	Recall the following form of Stirling's formula 	\cite[p. 224]{copson:1935}. Uniformly for $\sigma$ in any finite 	interval, as $|t|$ tends to $\infty$,
	\begin{equation}\label{E:stirling}
		|\Gamma(s)|\sim \sqrt{2\pi} e^{-\frac{\pi |t|}{2}} 		|t|^{\sigma-\frac{1}{2}}.
	\end{equation}
	Also, by \cite[p. 81]{titchmarsh:1951}, uniformly for 	$\sigma\ge\sigma_0$, there exists a constant 	$c=c(\sigma_0)>0$, such that
	\begin{equation}\label{E:titch}
		\zeta(s)=O(|t|^c),
	\end{equation}
	as $|t|\to\infty$. The identities \eqref{E:zeta1}--\eqref{E:zeta3} 	then imply that for any fixed $k$, uniformly for 	$3\ge\sigma\ge\sigma_0$, with the same $c$ as in 	\eqref{E:titch}
	\begin{equation}\label{E:bound3}
		\zeta((1+|k|)/2,s)=O(|t|^c),
	\end{equation}
	as $|t|\to\infty$. The estimates \eqref{E:bound1} and 	\eqref{E:bound2} now follow from \eqref{E:stirling} and 	\eqref{E:bound3}.
\end{proof}

\end{document}